\documentclass[12pt,reqno,twoside]{amsart}
\usepackage{amsmath}
\usepackage{amssymb}
\usepackage{epsfig}
\usepackage{color}

\numberwithin{equation}{section}

\title[
Dirichlet's Theorem and homogeneous flows]{
Dirichlet's theorem on diophantine approximation
and homogeneous flows
}

\author{Dmitry Kleinbock}

\address{Brandeis University, Waltham MA
02454-9110 {\tt kleinboc@brandeis.edu}}

\author{Barak Weiss}

\address{Ben Gurion University, Be'er Sheva, Israel 84105
{\tt barakw@math.bgu.ac.il}}

\newif\ifdraft\drafttrue

\draftfalse

\font\sb = cmbx8 scaled \magstep0

\font\sn = cmssi8 scaled \magstep0
\font\si = cmti8 scaled \magstep0

\long\def\comdima#1{\ifdraft{\si #1 }\else\ignorespaces\fi}

\long\def\combarak#1{\ifdraft{\sb #1 }\else\ignorespaces\fi}

\long\def\commargin#1{\ifdraft{\marginpar{\it #1}}\else\ignorespaces\fi}

\newcommand\name[1]{\label{#1}{\ifdraft{\sn [#1]}\else\ignorespaces\fi}}

\newcommand\eq[2]{{\ifdraft{\ \tt [#1]}\else\ignorespaces\fi}\begin{equation}\label{eq:
#1}{#2}\end{equation}}

\newcommand {\equ}[1]     {\eqref{eq: #1}}
\newcommand{\under}[2]{\underset{\text{#1}}{#2}}

\newcommand{\goth}[1]{{\mathfrak{#1}}}

\newcommand{\compose}{{\circ}}

\newcommand{\R}{{\mathbb{R}}}

\newcommand{\Z}{{\mathbb{Z}}}

\newcommand{\N}{{\mathbb{N}}}

\newcommand{\vf}{{\bf{f}}}
\newcommand{\vt}{{\bf{t}}}

\newcommand{\vv}{{\bf{v}}}
\newcommand{\vu}{{\bf{u}}}
\newcommand{\y}{{\bf{y}}}

\newcommand{\vr}{{\bf{r}}}

\newcommand{\vs}{{\bf{s}}}

\newcommand{\ve}{{\bf{e}}}

\newcommand{\vw}{{\bf{w}}}

\newcommand{\Id}{{\operatorname{Id}}}

\newcommand{\GL}{\operatorname{GL}}

\newcommand{\SL}{\operatorname{SL}}

\newcommand{\ggm}{G/\Gamma}

\newcommand{\diag}{{\operatorname{diag}}}
\newcommand{\vol}{{\operatorname{vol}}}

\newcommand{\spa}{{\rm span}}

\newcommand{\TT}{{\mathcal T}}

\newcommand{\df}{{\, \stackrel{\mathrm{def}}{=}\, }}

\newcommand{\s}{{\bf y}}

\newcommand{\x}{{\bf x}}

\newcommand{\vp}{{\bf p}}

\newcommand{\vm}{{\bf m}}
\newcommand{\vn}{{\bf n}}

\newcommand{\til}{\widetilde}

\newcommand{\supp}{{\rm supp}}

\newcommand{\vq}{{\mathbf{q}}}

\newcommand{\T}{{\mathbf{t}}}

\newcommand{\sm}{\smallsetminus}

\newcommand{\vre}{\varepsilon}

\newcommand{\nz}{\smallsetminus\{0\}}

\font\sb = cmbx8 scaled \magstep0

\newcommand\cag{$(C,\alpha)$-good}
\newcommand\hd{Hausdorff dimension}

\newcommand\ba{badly approximable}

\newcommand\da{diophantine approximation}

\newcommand\di{diophantine}
\newcommand\dc{Dirichlet constant}

\newcommand\hs{homogeneous space}

\newcommand{\fa}{{\frak a}}

\newcommand{\fr}{{\mathcal R}}

\newcommand{\DI}{{\mathrm{DI}}}

\newcommand\dt{Dirichlet's Theorem}

\newcommand\qn{quantitative nondivergence}

\newcommand\mr{M_{m,n}
}
\newcommand\amr{$Y\in M_{m,n}
$}

\newcommand {\ignore}[1]  {}

\newcommand\ssm{\smallsetminus}

\newtheorem{thm}{Theorem}[section]

\newtheorem{lem}[thm]{Lemma}

\newtheorem{prop}[thm]{Proposition}

\newtheorem{cor}[thm]{Corollary}

\date{September 2007}
\begin{document}

\ignore{

\begin{abstract}

\end{abstract}

}


\dedicatory{Dedicated to Gregory Margulis with admiration and respect}

\maketitle

\begin{abstract}
Given an $m \times n$ real matrix $Y$, an unbounded set $\mathcal{T}$
of parameters $\vt =\left( t_1, \ldots,
t_{m+n}\right)\in\R_+^{m+n}$ with $\sum_{i = 1}^m t_i =\sum_{j = 1}^{n} t_{m+j} $ and $0<\vre
\leq 1$, we say that Dirichlet's Theorem can be $\vre$-improved 
 for $Y$ along $\mathcal{T}$ if for every sufficiently large $\vt \in
\mathcal{T}$ there are nonzero $\vq \in \Z^n$ and $\vp \in \Z^m$ 
such that 
\[
\begin{cases}
|Y_i\vq - p_i| < \vre e^{-t_i}\,,\quad &i = 1,\dots,m
 \\  
\ \ |q_j| < \vre e^{t_{m+j}}\,,\quad &j = 1,\dots,n
\end{cases}
\]
(here $Y_1,\dots,Y_m$
are rows of $Y$).
We show that for any $\vre<1$ and 
any $\mathcal{T}$ `drifting away from walls', see \equ{drift}, 
 Dirichlet's Theorem cannot
be $\vre$-improved along $\mathcal{T}$ for Lebesgue almost every $Y$. 
In the case $m = 1$ we also show that for a large
class of measures $\mu$ (introduced in
\cite{friendly}) there is $\vre_0>0$ such that for 
any 
unbounded $\mathcal{T}$, any $\vre<\vre_0$, and for $\mu$-almost every $Y$,
Dirichlet's Theorem cannot
be $\vre$-improved along $\mathcal{T}$. These measures include 
natural measures on sufficiently regular smooth manifolds and fractals. 

Our results extend those of several authors beginning with the work of Davenport
and Schmidt done in late 1960s. The proofs rely on a translation of the problem
into a dynamical one regarding the action of a diagonal semigroup on
the space $\SL_{m+n}(\R)/\SL_{m+n}(\Z).$ 
\end{abstract}

\ignore{
\combarak{This was a rather difficult abstract to write given that
many definitions are needed. It is not completely precise because I
omitted 
the drifting away from walls assumption nor did I describe adequately
the `sufficiently regular manifolds and fractals'. Still I think it is precise
enough and informative enough, see if you want to change it. }
\comdima{My changes: removed the reference to Dirichlet since the statement with $<$ was
false;  mentioned the drift  and $m=1$ so that we don't make false claims, a slight
change of date in Davenport-Schmidt.}
}

\section{Introduction}

Let $m,n$ be positive integers, and denote by $\mr$ 
the space of $m\times n$ matrices with real entries. 
\dt\ (hereafter abbreviated by `DT') 
on simultaneous diophantine approximation 
states that for any \amr\ (viewed as a system
of $m$ linear forms in $n$ variables) 
and for any 
$t > 0$ there exist
$\vq = (q_1,\dots,q_n) \in \Z^n\nz$ and 
$\vp = (p_1,\dots,p_m)\in \Z^m$ satisfying the following system of inequalities:
\eq{dt}{
\|Y\vq - \vp\|
< e^{-t/m}
  \ \ \ \mathrm{and}  
\ \ \|\vq\|
\le e^{t/n}
\,.}
Here and hereafter, unless otherwise specified, $\|\cdot\|$ stands for the norm on $\R^k$ given by 
$\|\x\| = \max_{1\le i \le k}|x_i|$.  See \cite{Schmidt:book}
for a discussion of two ways of proving this theorem, due to
Dirichlet and Minkowski respectively.

Given $Y$ as above and positive $\vre < 1$, 
we will say that 
DT {\sl can be $\vre$-improved\/} for $Y$, 
and write $Y\in\DI_\vre(m,n)$, or  $Y\in\DI_\vre$ when the
dimensionality is clear from the context, 
if 
for every sufficiently
large $t$  one can find
$\vq  \in \Z^n\nz$ and 
$\vp \in \Z^m$ with
\eq{di}{
\|Y\vq - \vp\|
 < \vre e^{-t/m} 
  \ \ \ \mathrm{and}  
\ \ \|\vq\|
 < \vre  e^{t/n} 
\,,}
that is, satisfy \equ{dt} with the right hand side terms 
 multiplied by $\vre$ (for convenience we will also replace $\le$ in the second inequality by $<$).
\ignore{ \comdima{I am not sure if it makes sense to introduce this new notation,
let me know what you think}
Further, to measure the extent to which
DT can be improved for a given  $Y$, let us define the {\it \dc\/}
$\vre(Y)$ of $Y$ to be the infimum of $\vre$ for which $Y\in\DI_\vre$.
\begin{equation}
\label{eq: dc}
\vre(Y) = \inf\left\{ \vre:
\begin{aligned}  \text{ the system }\equ{di} \text{  has nontrivial integer}\\
\text{solutions for every
large enough } T\quad
\end{aligned}\right\}
\,.
\end{equation}
Thus DT implies that $\vre(Y) \le 1$ for all $Y$.} Also note that 
$Y$ is called {\sl singular\/} if 
$Y\in\DI_\vre$
for any $\vre > 0$ 
(in other words, DT can be `infinitely 
improved' for $Y$). 

The two papers \cite{Davenport-Schmidt, Davenport-Schmidt2}
by H.\ Davenport and W.\ Schmidt give a few basic results concerning 
the properties defined above. For example, the following is proved
there:

\begin{thm}[\cite{Davenport-Schmidt2}]\name{thm: ds1}
 For any\footnote{Even
though the results of 
 \cite{Davenport-Schmidt2} are stated for matrices with one row or
one column, i.e.\ for a vector or a single linear form,  the proofs
can be 
generalized to the setting of systems of linear forms.} 
$m,n\in\N$ and any $\,\vre < 1$, 
the sets  $\DI_\vre(m,n)$ 
have Lebesgue measure zero.
\end{thm}
In other words, $\lambda$-generic systems of linear forms do not allow 
any improvement to DT ($\lambda$ will denote Lebesgue measure throughout the paper).

Another question raised by  Davenport and Schmidt concerns
the possibility of improving DT for matrices with 
some functional
relationship between entries. Specifically they considered row
matrices
  \eq{ds}{\vf(x) = \begin{pmatrix}x & x^2\end{pmatrix}\in M_{1,2}\,,}
and proved 
\begin{thm}[\cite{Davenport-Schmidt}]\name{thm: ds2}
For any $\,\vre < 4^{-1/3}$, the set of $x\in \R$ for which 
$\vf(x)\in\DI_\vre(1,2)$ has zero Lebesgue measure. 
\end{thm}
In other words, generic matrices of the form \equ{ds} do not allow a  
sufficiently drastic improvement to DT.
This result was subsequently extended by R.\ Baker, Y.\ Bugeaud, and others.
Namely, for some other smooth submanifolds of  $\R^n$ they exhibited
constants $\vre_0$ such that almost no points on these submanifolds
(viewed as row or column matrices) are in $\DI_\vre
$ for $\vre < \vre_0$. 
We will discuss the history in more detail in
\S \ref{const}.
\medskip

In the present paper we significantly generalize Theorems \ref{thm: ds1} and
\ref{thm: ds2}
 by using a homogeneous dynamics approach and following
a theme developed in the paper \cite{sing} which dealt with
infinite improvement to DT, that is, with singular systems.  
Namely, 
we study the subject of 
improvement of the multiplicative version
of DT. In what follows we fix $m,n\in\N$ and let $k = m+n$.
Let us denote by $\fa^+
$
the set of $k$-tuples $\vt = (t_1,\dots,t_{k})\in \R^{k}$
such that
\eq{sumequal}{
t_1,\dots,t_{k} > 0
\quad \mathrm{and}\quad 
\sum_{i = 1}^m t_i =\sum_{j = 1}^{n} t_{m+j} \,.
} 
It is not hard to see that both Dirichlet's and Minkowski's proofs of 
DT easily yield the following statement:

\begin{thm}\name{thm: mdt} For any system  of $m$ linear forms $Y_1,\dots,Y_m$
(rows of \amr)
in $n$ variables
and for any $\,\vt \in\fa^+
$
there exist solutions $\vq  = (q_1,\dots,q_n)\in \Z^n\nz$ and 
$\vp = (p_1,\dots,p_m) \in \Z^m$ 
of
\eq{mdt}{
\begin{cases}
|Y_i\vq - p_i| < e^{-t_i}\,,\quad &i = 1,\dots,m
 \\  
\ \ |q_j| \le e^{t_{m+j}}\,,\quad &j = 1,\dots,n
\,.
\end{cases}}
\end{thm}
Now, given an unbounded subset $\mathcal{T}$ of $\fa^+
$
and positive $\vre < 1$, say that
{\sl DT can be $\vre$-improved for $Y$  along\/} $\mathcal{T}$, 
or $Y\in\DI_\vre(\mathcal{T})$,  
 if there is $t_0$ such that for every  $\vt =
(t_1,\dots,t_{k})\in\mathcal{T}$ with $\|\vt\| >t_0$, 
the  inequalities 
\eq{mdtw}{
\begin{cases}
|Y_i\vq - p_i| < \vre e^{-t_i}\,,\quad &i = 1,\dots,m
 \\  
\ \ |q_j| < \vre e^{t_{m+j}}\,,\quad &j = 1,\dots,n
\,.
\end{cases}
}
i.e., \equ{mdt} with the right hand side terms 
 multiplied by $\vre$,
 have nontrivial integer solutions. Clearly
$
\DI_\vre
 =  \DI_\vre(\fr
)
$
where \eq{def r}{\fr
 \df \left\{\left(\tfrac t m,\dots,\tfrac t m,\tfrac t n,\dots,
\tfrac t n\right)
: t > 0\right\}} is 
the `central ray' in $\fa^+
$. Also say that 
$Y$ is  
{\sl  singular along\/} $\mathcal{T}$ 
if it belongs to $\DI_\vre(\mathcal{T})$
for every positive $\vre$. The latter definition was introduced  in 
\cite{sing} for $\mathcal{T}$ contained in an arbitrary ray in  $\fa^+
$
emanating from the origin (the setup of {\em \da\ with weights\/}) in
the special case $n = 1$.

To state our first
result we need one more definition. Denote 
$$
\lfloor \vt\rfloor \df \min_{i = 1,\dots,k}t_i\,,
$$ 
and say that 
$\mathcal{T}\subset\fa^+$
{\sl drifts away from walls\/} if \commargin{definition changed}
\eq{drift}{\forall \,s > 0\quad\exists\,
\vt\in\mathcal{T}  \text{ with }
\lfloor \vt\rfloor > s\,.}
In other words, the distance
from $\vt\in \mathcal{T}$ to the boundary of  $\fa^+$ 
is unbounded (hence, in partucular, $\mathcal{T}$ itself is unbounded).

Using Lebesgue's  Density Theorem and an elementary argument  which
can be found e.g.\ in  \cite[Chapter V, \S 7]{Cassels} 
and dates back to  Khintchine, one can show that 
for any $m,n$ and $\mathcal{T}\subset\fa^+
$ drifting away from walls,  
$\DI_\vre(\mathcal{T})$ has Lebesgue 
measure zero as long as $\vre < 1/2
$. 
However proving  the 
same for any $\vre < 1$ 
is more difficult, even in the case 
$\mathcal{T} = \fr
$ considered by Davenport and Schmidt. In \cite{Davenport-Schmidt2}
they derived this fact from the density of generic trajectories 
of certain one-parameter subgroups on the space
$\ggm$, where 
\eq{defn ggm}{G = \SL_{k}(\R)\text{ and }\Gamma = \SL_{k}(\Z)\,.} 
In \S \ref{unif} we show how the fact that 
$\DI_\vre(\mathcal{T})$ has Lebesgue 
measure zero 
for any $\vre < 1$ and any unbounded
$\mathcal{T}\subset\fr
$ 
follows from mixing of the $G$-action on $\ggm$.  
We also explain why mixing is not enough to obtain such a result
for an arbitrary $\mathcal{T}\subset\fa^+
$ drifting away from walls,  
even contained in a single `non-central' ray, 
and 
prove the following multiplicative analogue of Theorem~\ref{thm: ds1}:

\begin{thm}\name{thm: dsw} For any\,
$\mathcal{T}\subset\fa^+
$ drifting away from walls 
and any $\vre < 1$,
the set  $\DI_\vre(\mathcal{T})$ 
has Lebesgue measure zero.
\end{thm}

In other words, given any 
$\mathcal{T}$ as above, DT in its generalized form (Theorem \ref{thm: mdt})
cannot be improved for generic
systems of linear forms. 
 Theorem \ref{thm: dsw} is derived from 
the equidistribution of translates of certain measures on $\ggm$ (Theorem 
\ref{thm: unif distr}). 
Our proof of this theorem, described in \S \ref{dynamical}, is 
a modification of arguments used in analyzing epimorphic groups  \cite{epi, with Nimish}.
It relies on S.\,G. Dani's classification of measures invariant under horospherical
subgroups, and the `linearization method' developed by Dani, G.\,A.\
Margulis, M. Ratner,  J. Smillie, 
N.\ Shah and others. See also \cite{KM new} for an alternative approach
suggested to the authors by Margulis.

\medskip

Next, building on the approach of \cite{friendly, sing}, we generalize 
Theorem~\ref{thm: ds2}, namely, consider measures other than
Lebesgue. Here we will restrict ourselves to 
measures 
$\mu$ 
on $\R^n\cong M_{1,n}$ and study \di\ properties of
$\mu$-almost all 
$\y\in
\R^n$  interpreted as row vectors (linear forms);
that is, we put $m = 1$ and $k = n+1$. 
The dual case of simultaneous approximation
 ($n = 1$, which was the setup
of \cite{friendly} and \cite{sing}) can be treated along the same lines. 

All measures on Euclidean spaces will be assumed to be Radon (locally finite regular
Borel). Generalizing the setup of 
 Theorem~\ref{thm: ds2}, we will consider measures  on $\R^n$  of the form 
 $\vf_*\nu$, where $\nu$ is a measure on $\R^d$  and
 $\vf$ a map from $\R^d$ to $\R^n$. Our assumptions on $\vf$ and $\nu$ 
rely on definitions  of:

\pagebreak

\begin{itemize}
\item 
 measures $\nu$ which are {\sl $D$-Federer\/} on open subsets $U$ of $ \R^d$,
and 
\item 
functions $f:U\to \R$ which are {\sl \cag\ on $U$ w.r.t.\/} $\nu$ 
\end{itemize} 
(here $D,C,\alpha$
are positive constants). 
%
We postpone the precise definitions, introduced in \cite{KM, friendly},
until \S \ref{quant}. 

Now given a measure $\nu$  on $\R^d$, an open $U\subset \R^d$ with $\nu(U) > 0$  and
 a map  $\vf:\R^d\to \R^n$,  say that 
a pair $(\vf,\nu)$ is

\begin{itemize}
\item {\sl
\cag\  
on $U$\/}  if 
any linear 
combination of
$1,f_1,\dots,f_n$  is
$(C,\alpha)$-good on
$U$ with respect to $ \nu$;

\item {\sl 
nonplanar
on $U$\/}
    if
for any  ball
$B\subset U$ centered in $\supp\,\nu$, 
the
restrictions of $ 1,f_1,\dots,f_{n}$ to $B\,\cap \,\supp\,\nu$
are linearly independent over $\R$;
in other words, if $\vf(B \,\cap\,
\supp\,\nu)$  is not contained in any proper affine subspace  of $\R^n$;

\item   {\sl 
$(C,\alpha)$-good\/} (resp., {\sl 
nonplanar\/}) if  for $\nu$-a.e.\ 
$\x$  there exists  a neighborhood $U$ of $\x$  
such that $\nu$
is
 $(C,\alpha)$-good (resp., nonplanar) on $U$.

\end{itemize}

Similarly, we will say that $\nu$   is {\sl 
$D$-Federer\/} if  for $\nu$-a.e.\ 
$\x\in\R^d$  there exists a neighborhood $U$ of $\x$  
such that $\nu$
is
 $D$-Federer 
on $U$.

In \S\ref{est} we  prove

\begin{thm}\name{thm: friendly general} For any\,
$d,n\in\N$ and   
$\,C,\alpha,D > 0$ there exists
$\vre_0 =  \vre_0(d,n,C,\alpha,D)$ with the following property.
 Let 
$\nu$ be a  
measure  on
$\R^d$   and  $\vf$ a continuous map from open 
$U\subset\R^d$ with $\nu(U) > 0$ to $\R^n$. 
Assume that $\nu$ is  $D$-Federer,  and $(\vf,\nu)$ is 
$(C,\alpha)$-good and
nonplanar.
Then for any $\vre < \vre_0$
\eq{result}{ \vf_*\nu\big(\DI_\vre(\mathcal{T})\big) = 0\ 
\text{ for any unbounded } 
\mathcal{T}\subset\fa^+\,.
}
\end{thm}

One can easily check, see \S\ref{polynom}, that 
$(\vf,\lambda)$
as in \equ{ds} satisfies the assumptions of the above theorem.
Moreover, Theorem \ref{thm: friendly general} generalizes all known
results 
on improvement of DT on manifolds, 
and is applicable to   many other examples of interesting measures, see \S\ref{fr}
for more examples.
The 
proof 
is based on so-called `quantitative nondivergence' technique, 
a generalization of non-divergence of unipotent flows on \hs s first established by Margulis in 
early 1970s. 

Note that  $\vre_0$ in the above theorem can be explicitly computed in terms of 
$d,n,C,\alpha,D$, although the value that can be obtained  seems to be far from optimal.
In particular they are not as good as the values of $\vre_0$ 
obtained in \cite{Davenport-Schmidt, DRV1, Bugeaud}.
On the other hand the setup of Theorem \ref{thm: friendly general} is much more 
general than that of the aforementioned papers, as illustrated by examples 
considered in \S\ref{const}.

We remark  that the above theorem does not require
to impose an additional condition \equ{drift} on $\TT$, unlike Theorem \ref{thm: dsw}
which is in general false if  \equ{drift} is not assumed. Counterexamples are discussed in 
 \S\ref{nodrift}.

\medskip

{\bf Acknowledgements:}
This research was supported by BSF grant
2004149, ISF grant 584/04 and NSF grant DMS-0239463.
Thanks are due  to
Roger
Baker, Yann Bugeaud and Gregory Margulis for helpful 
discussions, and to the referee for useful remarks resulting in strengthening some of 
the results of the paper.

\section{Theorem \ref{thm: dsw} and equidistribution}
\name{dynamical}
\subsection{The correspondence}\name{corr}
Let $G$ and $\Gamma$ be as in \equ{defn ggm}, and denote by $\pi: G
\to \ggm$
the quotient map.
    $G$ acts on $\ggm$ by left translations via the rule $g \pi(h) =
\pi(gh)$, $g,h\in G$.
Define
$$
\tau(Y) \stackrel{\mathrm{def}}{=} \left(
\begin{array}{ccccc}
I_{m} & Y \\ 0 & I_{n} 
\end{array}
\right), \ \ \ \bar{\tau} \stackrel{\mathrm{def}}{=} \pi \compose \tau\,,
$$
where $I_{\ell}$ stands for the $\ell\times \ell$ identity matrix.
\ignore{Then, given $\vr,\vs$ as in \equ{setting r}, \equ{setting s}, 
 consider the  one-parameter subgroup
$\{g_t^{(\vr,\vs)}\}$ of $G$
given by
\begin{equation*}
\label{eq: new defn g_t}
g_t^{(\vr,\vs)} \df \diag(e^{r_1t}, \ldots,
e^{r_{m} t}, e^{-s_1t}, \ldots,
e^{-s_{n} t})\,.
\end{equation*}
}
Since $\Gamma$ is the stabilizer of $\Z^{k}$ under the action of $G$ on
the set of lattices in $\R^{k}$, $\ggm$ can be identified with 
$G\Z^{k}$,
that is, with the set of all unimodular lattices in $\R^{k}$.
To highlight the relevance of the objects defined above to the \di\ problems considered in the introduction,
note that 
\eq{expl tau}{
 \bar{\tau}(Y) = \left\{\begin{pmatrix} Y\vq - \vp\\\vq\end{pmatrix} : \vp\in\Z^m,\  \vq\in\Z^n\right\}\,.
}
Now for $\varepsilon>0$ let
\eq{keps}{
\begin{split}
K_{{\varepsilon}} &\stackrel{\mathrm{def}}{=} \pi\big(\big\{ g \in G : \Vert g
\vv \Vert \geq {\varepsilon} \quad \forall\, \vv \in
\Z^{k} \sm \{ 0 \}\big\}\big),
\end{split}}
i.e., $K_{ \varepsilon}$ is the collection of all unimodular
lattices in $\R^{k}$ which contain no nonzero vector 
of norm smaller than
$\varepsilon$.
By Mahler's compactness criterion (see e.g.\  \cite[Chapter
10]{Rag}), each $K_{{\varepsilon}}$ is compact,
and for each compact $K \subset \ggm$  there is $\vre>0$
such that $K \subset K_{ \vre}$. Note also that $K_\vre$
is empty
if $\vre > 1$ by Minkowski's Lemma, see \cite{Schmidt:book}, and has nonempty interior if $\vre <1$.

Now to any $\vt\in\fa^+
$ 
let us associate the diagonal
matrix $$g_\vt \df \diag(e^{t_1}, \ldots,
e^{ t_m}, e^{-t_{m+1}}, \ldots,
e^{-t_{k}})\in G\,.$$
Then, using \equ{expl tau}, it is straightforward to see  that the system \equ{mdtw} has a nonzero
integer solution if and only if 
$g_\vt \bar{\tau}(Y)\notin K_{\vre}$.
 We therefore arrive at

\begin{prop}
\name{prop: dynamical interpretation}
For \amr, $0 < \vre < 1$ and unbounded $\mathcal{T}\subset\fa^+
$, 
one has $Y\in\DI_\vre(\mathcal{T})$ if and only if $g_\vt \bar{\tau}(Y)$ is outside of $K_{\vre}$ 
 for all $\vt\in\mathcal{T}$
with large enough norm. Equivalently,
\eq{liminf}
{\DI_\vre(\mathcal{T}) = \bigcup_{t > 0}\quad\bigcap_{\vt\in\mathcal{T} ,\, \|\vt\| > t}
\{Y : 
g_\vt \bar{\tau}(Y)\notin K_{\vre}\}\,.}
\commargin{replaced $t_0$ by $t$}
\end{prop}

In particular, $Y$ is singular along $\mathcal{T}$ if and only if  
the trajectory $\{g_\vt \bar{\tau}(Y): \vt\in\mathcal{T}\}$ is divergent 
(i.e.\ eventually leaves $K_{\vre}$ for any $\vre > 0$). The latter observation was first
made in \cite[Proposition 2.12]{Dani} for $\mathcal{T} = \fr
$,
and then in \cite[Theorem 7.4]{dima duke} for an arbitrary ray in $\fa^+
$.



\subsection{Lebesgue measure and uniform distribution}\name{unif}
We recall that one of the goals of this paper is to show that,
whenever $\vre < 1$ and $\TT\subset\fa^+$ drifts away from walls, $\lambda$-a.e.\
\amr\ does not belong to  
$\DI_\vre(\mathcal{T})$ 
(Theorem \ref{thm: dsw}). 
\commargin{Shortened this paragraph and killed (2.4)}
The latter, in view of the above proposition, amounts to  showing that 
the set $\{\vt\in\mathcal{T} : g_\vt \bar{\tau}(Y) \in K_\vre \}$ is unbounded
 for $\lambda$-a.e.\ $Y\in\mr$. 
In this section we will prove a
stronger statement. Here and hereafter 
`$\vol$'  stands for  the $G$-invariant probability measure  on $\ggm$.

\begin{thm}
\name{thm: unif distr} 
\commargin{Changed $\vre$ to $\delta$ not to mix with $K_\vre$}
Let a continuous compactly supported function $\varphi$ on $\ggm$ 
and bounded $B\subset \mr$ with positive 
Lebesgue measure be given. 
Then for any compact subset 
$L$ of $G/\Gamma$ and any $\,\delta>0$ there exists $T > 0$ such that 
\eq{unifdist}{
\left|\frac1{\lambda(B)}\int_{B}
\varphi\big(g_{\vt}\tau(Y){z}\big)
\, d\lambda(Y) - 
\int_{\ggm} \varphi\, d\vol\right| < \delta\,
}
for all ${z}\in L$ and $\vt\in\fa^+, \, \lfloor \vt\rfloor\ge T$.
\end{thm}

Let us denote by 
$\lambda_{{z},B}$ the pushforward of $\frac1{\lambda(B)}\lambda|_B$ to 
$\ggm$ by the map
$Y\mapsto \tau(Y){z}$. 
That is, for a Borel subset $A$ of $\ggm$, let
$$\lambda_{{z},B}(A)  \df \frac{\lambda \big(\left\{Y \in B: \tau(Y){z} \in A
\right\}\big)}{\lambda(B)}.$$ 
The above theorem asserts that 
$g_{\vt}$-translates of $\lambda_{{z},B}$ weak-$*$ converge to 
$\vol$ as $\lfloor
 \vt\rfloor
\to\infty$, and the convergence  is uniform in ${z}$ when the latter
is restricted to a compact subset of $\ggm$.  
By taking ${z}$ to be the standard lattice $\Z^{k}\in \ggm$
and approximating $K_\vre$  
(which has nonempty interior and boundary of measure zero for
any $\vre < 1$) 
by continuous functions on $\ggm$, 
from Theorem \ref{thm: unif distr} one obtains that 
$$
\lambda\big(
\{Y\in B : g_{\vt} \bar{\tau}(Y)
\in K_\vre\} \big) \to  \lambda(B)\vol(K_\vre)\quad\text{ as }\lfloor
 \vt\rfloor
\to\infty
\,.
$$
This immediately rules out the existence of $\mathcal T\subset\fa^+$
 drifting away from walls and 
$B\subset\mr$ of positive Lebesgue measure
such that the set $\{\vt\in\mathcal{T} : g_\vt \bar{\tau}(Y) \in K_\vre \}$ is bounded for  any $Y\in B$. 
Thus 
Theorem \ref{thm: dsw}
follows from  Theorem \ref{thm: unif distr}.
\ignore{
 for any subset
$B
\subset \mr$ of positive measure there exists $t_0$ such that 
$$t \geq t_0 \ \ \implies \ \ g_{t}^{(\vr,\vs)} \bar{\tau}(B) \cap K_\vre
\neq
\varnothing.$$
This readily implies Theorem \ref{thm: dsw dynamical}, and hence Theorem
\ref{thm: dsw}.}

\medskip

Note that the conclusion of  Theorem \ref{thm: unif distr} is not new for $\mathcal{T}
$ being a subset of the `central ray' $\fr$ defined in \equ{def r}.
Indeed, one can immediately see that the group
 \eq{defn u}{H\df\tau(\mr)}
 is 
{\em expanding horospherical\/} (see e.g.\ \cite[Chapter 1]{handbook} for the definition)
with respect to $g_{\vt}$ for any $\vt\in\fr$.
In this case  Theorem \ref{thm: unif distr}
coincides with \cite[Proposition 2.2.1]{KlMa}, where it was
deduced from mixing of the $G$-action on
$\ggm$, using an argument 
dating back to the Ph.D.\ Thesis of Margulis \cite{Mar:Anosov}.
On the other hand, $H$ is strictly contained
in the expanding horospherical subgroup relative to $g_{\vt}$
for any $\vt\in\fa^+\ssm\fr $,
thus the aforementioned argument does not prove Theorem \ref{thm: unif distr}.
In the remainder part of this section 
we show how to bypass this difficulty. See also \cite{KM new} for an alternative proof.

\subsection{Expanding vectors along cones}
\name{drift}
In this section we will discuss an important representation-theoretic property of
the 
pair $(\fa^+, H)$.

\begin{lem}
\name{lem: main computation}
Let $\rho: G \to \GL(V)$ be a representation (of algebraic groups)
on a finite-dimensional normed vector space $V$ 
without nonzero fixed vectors,
and let 
$$V^H = \{v \in V:  \rho(h)v=v\ \forall h \in H\}.$$
Then there are positive $c, c_0$ such that for any $v
\in V^H$ and $\vt\in\fa^+$ one has 
$$
\|\rho(g_{\vt})v\| \ge c_0e^{c\lfloor \vt\rfloor}\|v\|\,.
$$
\end{lem}

\begin{proof}
Let us denote by $A$ the group of positive diagonal matrices in $G$,
let $\goth{a}$ be its Lie algebra, and let $\goth{h}$ be the Lie
algebra of $H$. Since $A$ normalizes $H$, $V^H$
is a $\rho(A)$-invariant subspace, and we may write
$$V^H = \bigoplus_{\chi \in \Psi} V_{\chi},$$
where $\Psi$ is a finite set of weights (linear functionals on
$\goth{a}$) and
$$V_{\chi} = \{v \in V : \rho(\exp Y)v =
e^{\chi(Y)}v \ \forall\, Y \in \goth{a}\}$$
is nonzero for any $\chi \in \Psi$. There is no loss of generality in
assuming that $\| \cdot \|$ is the sup-norm with respect to a basis of
$\rho(A)$-eigenvectors; thus it suffices to show that for
any $\chi \in \Psi$, 
$\inf_{\vt\in\fa^+}\chi(\vt)/ \lfloor \vt\rfloor$ is positive.

Let 
$$\mathcal{I} = 
\{1, \ldots, m\} \times \{m+1, \ldots, k\}
,$$
and for each $(i,j) \in \mathcal{I}$ let $G_0 = G_0(i,j)$ be the Lie
subgroup of $G$ whose Lie algebra $\goth{g}_0$ is 
generated by $E_{ij}, E_{ji}, F_{ij}$ (here $E_{rs}$ stands
for the $k \times k$ matrix with $1$ in position $(r,s)$ and $0$
elsewhere, and $F_{rs}$ stands for $E_{rr}-E_{ss}$). Then $G_0$ is a
copy of $\SL_2( \R)$ contained in $G$, such that 
$\goth{h} \cap \goth{g}_0 = \spa \, E_{ij}$ and $\goth{a} \cap
\goth{g}_0 = \spa \, F_{ij}.$ Since all vectors in $V_{\chi}$ are
 fixed by $\rho\big(\tau(E_{ij})\big)$, by  the representation theory of
$\goth{sl}_2$ (see e.g. \cite{Serre}) we have 
$\chi(F_{ij}) \geq 0$, and 
$$\chi(F_{ij})=0 \ \ \Longleftrightarrow \rho(g)v=v \  \forall\, g \in G_0,
\ v \in V_{\chi} \,.$$
Since $V$ contains no nonzero vectors fixed by $\rho(G)$, and since the
group generated by $\{G_0(i,j) : (i,j) \in \mathcal{I}\}$ is equal to
$G$,
there is at least one $(i_0,j_0)$ for which $\chi(F_{i_0j_0})
>0$. Since any $Y \in \overline{\goth{a}^+}$ can be written as a linear
combination of $\left\{F_{ij}: (i,j) \in \mathcal{I} \right \}$ with
non-negative coefficients, and since $\vt_0=\vt - \lfloor \vt \rfloor
(F_{i_0j_0}) \in \overline{\goth{a}^+}$, we have: 
$$\chi(\vt) = \chi(\vt_0) + \chi(\lfloor \vt \rfloor F_{i_0j_0}) \geq  \lfloor \vt \rfloor
\chi(F_{i_0j_0})\,,$$
finishing the proof.
\end{proof}
To put this result in context, recall that 
a subgroup $L_1$ of an algebraic group $L_2$ is said to be {\em
epimorphic} in $L_2$ if for any representation $\rho: L_2 \to \GL(V)$, any vector
fixed by $\rho(L_1)$ is also fixed by $\rho(L_2)$. For example, in our
present notation, $AH$ is
epimorphic in $G$. The `cone lemma' (\cite[Lemma 1]{epi}) shows that if $TU$ is
epimorphic in $L$, where $U$ is unipotent, $T$ is diagonalizable and
normalizes $U$, and $L$ is generated by unipotents, then for any
$\rho: L \to \GL(V)$ without nonzero 
fixed vectors there is a nonempty open
cone $T^+$ in $T$ such that 
for any nonzero $v\in V^U$, $\rho(a)v \to \infty$ as $a \to \infty$ in $T^+$. The 
proof in \cite{epi} is non-constructive. The true meaning of Lemma
\ref{lem: main computation} is a precise determination of a cone 
which works for
all representations $\rho$ in the case $T = A$, $U = H$. 

\medskip

The next proposition is a consequence of  Lemma
\ref{lem: main computation}.

\begin{prop}
\name{prop: repn}
Let $V$, $\rho$, $c$ 
be as in Lemma
\ref{lem: main computation}, and 
let $B$ be a neighborhood of 
$0$ in $\mr$. Then there exists $b 
> 0$ such that 
for any $v \in V$
 and $\vt\in\fa^+$ 
one has 
\eq{goes to infinity}{
\sup_{Y \in B} \big\| \rho\big(g_{\vt}\tau(Y)\big)v \big\| \ge 
be^{c\lfloor \vt\rfloor}\|v\|\,.
}
\end{prop}
\begin{proof} 
Denote by $p$ the $\rho(A)$-equivariant projection $V\to V^H$.
By \cite[Lemma 5.1]{Shah Indian academy} there exists $c_1 > 0$ (dependent on
$\rho$, $B$ and the choice of the norm) such that for any $v \in V$,
\eq{nimish}{\sup_{Y \in B} \big\| p\circ\rho\big(\tau(Y)\big)v \big\|\ge c_1\|v\|\,.}
Also choose $c_2 > 0$ such that $\|v\| \ge c_2\| p(v) \|$ for all $v \in V$.
Then for any $Y \in B$, $v \in V$ and $\vt\in\fa^+$ one can write
\begin{equation*}
\begin{aligned}
 \big\| \rho\big(g_{\vt}\tau(Y)\big)v \big\| \ge  &\ c_2 \big\| p\circ\rho\big(g_{\vt}\tau(Y)\big)v\big\| = 
 c_2 \big\|\rho(g_{\vt})\circ p\circ\rho\big(\tau(Y)\big)v \big\|\\  
\under{by Lemma \ref{lem: main computation}}\ge &\ c_0c_2 e^{c\lfloor \vt\rfloor}
\big\|  p\circ\rho\big(\tau(Y)\big)v \big\|\,,
\end{aligned}
\end{equation*}
hence
$$
\sup_{Y \in B} \|\rho\big(g_{\vt}\tau(Y)\big)v \| \ge c_0c_2 e^{c\lfloor \vt\rfloor}
\sup_{Y \in B}\big\| p\circ\rho\big(\tau(Y)\big)v \big\| \under{by \equ{nimish}}\ge  
c_0c_1c_2 e^{c\lfloor \vt\rfloor}\|v\|\,.
$$
\end{proof}

\subsection{Recurrence to compact sets}\name{recur} In order to 
establish the equidistribution of $g_{\vt}$-translates of 
$\lambda_{{z},B}$ as $\lfloor
 \vt\rfloor
\to\infty$, 
\ignore{\commargin{Perhaps some historical comment 
is appropriate here as to how nondivergence is an ingredient for 
topological Ratner stuff \combarak{I decided against it, we have too
much background and history as it is.}}}
one needs to at least show the existence
of one limit point (which is not guaranteed apriori since $\ggm$ is not compact).
In other words, there must exist a compact subset $K$ of $\ggm$ 
such that $\lambda_{{z},B}\big(g_{\vt}^{-1}(K)\big)$ is big enough whenever $\lfloor
 \vt\rfloor$ is large. 
We show in this section how to construct such a compact set
using Proposition \ref{prop: repn} and a theorem of Dani and Margulis.

Denote by $\goth{g}$ the Lie algebra of $G$, let 
$V \df \bigoplus_{j=1}^{\dim(G)-1} \bigwedge^j \goth{g}$, and let $\rho: G
\to \GL(V)$ be the representation obtained by acting on $V$ via the
adjoint representation and its exterior powers. 
Note that $V$ has no
nonzero $G$-fixed vectors since $G$ is simple. For any proper
connected Lie subgroup $W$ of $G$ we will denote by $\mathbf{p}_W$ an
associated vector in $V$.

We have the following result of Dani and Margulis, see  
\cite[Thm.\ 2.2]{Shah Indian academy} for a more general statement:
\ignore{\commargin{Is this the best reference? \combarak{this seems closest to
what we need}}}
\begin{prop}
\name{prop: DM nondivergence}
\ignore{\commargin{Also, maybe we should make some more comments,
as far as I understand this theorem has been stated in a very different disguise,
at least definitely without our notation $g_\vt$ and
$\tau$... \combarak{true but a dedicated reader should be able to
unravel it.}}}
Let $G$, $\Gamma$ and $\pi: G \to G/\Gamma$ be as above. Then there exist
finitely many closed subgroups
$W_1, \ldots, W_{\ell}$ of $G$ such that $\pi(W_i)$ is compact and
$\rho(\Gamma) \mathbf{p}_{W_i}$ is discrete for each $i \in \{1,
\ldots, \ell\}$, and the following holds: for any positive $\alpha, \vre$
there is a compact $K \subset \ggm$ such that for any $g \in G$, $\vt
\in \goth{a}$ and bounded convex open $B \subset M_{m,n}$, one of the following
is satisfied:
\begin{enumerate}
\item
There is $\gamma \in \Gamma$ and $i \in \{1, \ldots, \ell\}$ such that 
$${\sup_{Y \in B} \|\rho\big(g_{\vt} \tau(Y)g\gamma\big) \mathbf{p}_{W_i} \| < \alpha\,.
}
$$
\item $(g_{\vt})_*\lambda_{\pi(g),B}(K) \geq 1-\vre$.
\end{enumerate}
\end{prop}

\begin{cor}
\name{cor: recurrence}
For any compact subset 
$L$ of $G/\Gamma$ and any $\vre>0$  there exists a compact 
$K\subset \ggm$  with the following property: for any  bounded convex open
neighborhood $B$ of $0$
in $\mr$
there exists $T > 0$ such that 
$(g_{\vt})_*\lambda_{{z},B}(K) \geq 1-\vre$ whenever ${z}\in L$ and 
$\lfloor\vt\rfloor \ge T$.
\end{cor}

\begin{proof} Let $W_1, \ldots, W_{\ell}$ be as in 
Proposition \ref{prop: DM nondivergence}, and for any compact $L\subset \ggm$ 
consider
$$
\delta(L) \df \inf_{\pi(g)\in L,\ \gamma\in\Gamma,\ i = 1,\dots,\ell}\|\rho(g\gamma) 
\mathbf{p}_{W_i}\|\,.
$$
It is positive since $L$ is compact and 
$\rho(\Gamma) \mathbf{p}_{W_i}$ is discrete.
Proposition \ref{prop: repn} then 
implies that for any neighborhood $B$ of $0$
in $\mr$ there exist  constants $b,c$ such that for any 
$g \in \pi^{-1}(L)$, $\gamma\in\Gamma$, $i = 1,\dots,\ell$ and $\vt\in\fa^+$,
one has
\eq{big}{\sup_{Y \in B} \|\rho\big(g_{\vt} \tau(Y)\big) \rho(g\gamma)\mathbf{p}_{W_i} 
\|  \ge 
be^{c\lfloor\vt\rfloor}\delta(L)\,.}
Now take an arbitrary $\vre > 0$ and   $\alpha = 1$,
and choose $K$ according to Proposition \ref{prop: DM nondivergence}.
Then it follows from \equ{big} that for any $B$ there exists $T$ such that
whenever $\lfloor\vt\rfloor \ge T$ and $\pi(g) \in L$, the second
alternative of Proposition \ref{prop: DM nondivergence} 
must hold.
\end{proof}

\subsection{The linearization method}
As remarked in the introduction, our 
proof of Theorem \ref{thm: unif distr} relies on the work of many
mathematicians. Although we do not require Ratner's results on the
classification of measures invariant under unipotent flows (the earlier
results of Dani on horospherical subgroups are sufficient for us), we do use
the linearization method developed by many authors following Ratner's work.  
These results are described in detail in 
\cite{handbook}. Since our argument will be
very close to arguments in 
\cite{Shah Indian academy, with Nimish} we will rely on the notation
and results as stated in \cite{with Nimish}, where additional references to
the literature may be found. 

Let $\mathcal{H}$ be the set of all closed connected subgroups $W$
of $G$ such that $W \cap \Gamma$ is a lattice in $W$, and the subgroup of
$W$ generated by its one-parameter unipotent subgroups acts
ergodically on $W /(W \cap \Gamma)$. This is a countable collection. For
any $W \in \mathcal{H}$, we define 
\[
\begin{split}
 N(W, H) & = \{g \in G: Hg \subset gW\} \\
S(W,H) & = \bigcup \big\{N(W',H): W' \in \mathcal{H}, W' \subset W, \dim
W' < \dim W \big \} \\
N^*(W,H) & = N(W, H) \sm S(W, H).
\end{split}
\]
%
%

Recall that the subgroup $H$ of $G$ is horospherical. Dani \cite{Dani:
  horospherical} classified all the measures on $\ggm$ invariant under the
  $H$-action. 
The following is a consequence of Dani's 
classification  
and ergodic
decomposition:
\begin{prop}
\name{prop: Ratner}
Let $\mu$ be a finite $H$-invariant measure on $\ggm$ which is not
equal to $\mathrm{vol}$.
Then $\mu\big(\pi(N(W,H))\big)>0$ for some $W \in \mathcal{H}$ which is a
proper subgroup of $G$. 
\end{prop}

For $W \in \mathcal{H}$ let $V_W$ be the span of $\rho(N(W,H))
\mathbf{p}_W$ in $V$ and let $N^1_G(W) = \{g \in G: \rho(g)\mathbf{p}_W
= \mathbf{p}_W\}.$ Then the orbit $\rho(\Gamma) \mathbf{p}_W$ is
discrete in $V$. 

The next proposition uses the representation $\rho$ defined
in \S \ref{recur} to detect orbits which stay close to $\pi\big(N(W,H)\big)$
for some $W \in \mathcal{H}$. The idea has a long history and is
used in a similar context by Dani and Margulis in \cite{DM}. 

\begin{prop}
\name{prop: linearization} \cite[Theorem 4.1]{Shah Indian academy}
Given $W \in \mathcal{H}$ and $\vre>0$, 
for any compact $C \subset
\pi(N^*(W,H))$ there exists a compact $D \subset V_W$ with the
following property: for any neighborhood $D'$ of $D$ in $V_W$ there
exists a neighborhood $C'$ of $C$ in $\ggm$ such that for any $g \in
G$, any $\vt \in \goth{a}^+$ and any bounded convex open $B \subset
M_{m,n}$, one of the following holds:
\begin{enumerate}
\item There is $\gamma \in \Gamma$ such that 
$\rho(g_{\vt}\tau(B)g\gamma) \mathbf{p}_W \subset D'.$
\item $(g_{\vt})_*\lambda_{\pi(g),B}(C') 
< \vre.$
\end{enumerate}
\end{prop}

\subsection{Proof of Theorem \ref{thm: unif distr}}
Since finite linear combinations of indicator functions of 
balls are dense in $L^1(\mr,\lambda)$, there is no loss of generality in 
assuming that $B$ is an open ball in $\mr$. Take a sequence 
of points
${z}_n\in L$ and a sequence 
$\vt_n\in\fa^+$ drifting away from walls. 
It follows from Corollary \ref{cor: recurrence}
that  the sequence of translated measures  
$(g_{\vt_n})_*\lambda_{{z}_n,B}$  is weak-$*$ precompact, that is, 
along a subsequence we have $g_{\vt_n} \lambda_{{z}_n,B}$ to $\mu$, where 
$\mu$ is a Borel probability measure on $\ggm$. Our
goal is thus to show that $\mu = \mathrm{vol}$.

The hypothesis about drifting away from walls implies that for any $h
\in H$ we have $g_{\vt_n}^{-1} h g_{\vt_n} \to e,$ where $e$ is the
identity element in $H$. A simple computation (see \cite[Claim
3.2]{with Nimish}) shows that $\mu$ is $H$-invariant. By Proposition
\ref{prop: Ratner}, if $\mu \neq \mathrm{vol}$, there is a proper
subgroup $W \in \mathcal{H}$ such that $\mu\big(\pi(N(W,H))\big)>0.$
Making $W$ smaller if necessary we can assume that
$\mu\big(\pi(N^*(W,H))\big)>0.$ Let $C \subset \pi\big(N^*(W, H)\big)$
be compact with $\mu(C)>0,$ and put $\vre \df \mu(C)/2$.
Let $\til L \subset G$ be a compact subset such that $\pi(\til L) =
L$, and let $g_n \in \pi^{-1}({z}_n) \cap \til L$.  Applying Proposition
\ref{prop: linearization}, we find that there is a compact $D \subset
V_W$ such that the following holds. 
For each $n$, let ${D}_{n+1}
\subset {D}_n$ be a compact neighborhood of $D$ in $V_W$ such that
$\bigcap_n {D}_n =D$. Then there is an open neighborhood ${C}_n$ of
$C$ in $\ggm$ such that one of the following holds: 
\begin{enumerate} 
\item 
there is 
$v_n \in \rho(g_n\Gamma) \mathbf{p}_W$ such that $\rho(g_{\vt_n}
\tau(B)) v_n \subset {D}_n$; 
\item 
$(g_{\vt_n})_*\lambda_{{z}_n,B}({C}_n)
< \vre.$ 
\end{enumerate} 
Since $g_{\vt_n} \lambda_{{z}_n,B} \to \mu,$ and the sets ${C}_n$ are
neighborhoods of $C$, we find that $\lambda_{{z}_n,B} (g_{\vt_n}^{-1}
{C}_n) > \mu(C)/2 = \vre$ for all sufficiently large $n$,
so condition (2) above does not hold.  Therefore
$$\rho\big(g_{\vt_n} \tau(B)\big) v_n \subset {D}_n \subset {D}_1\,,$$
a bounded subset of $V_W$. 
On the other hand, since $\til L$ is
compact and $\rho(\Gamma)\mathbf{p}_W$ is discrete, we have 
$$\inf_{n} \|v_n\|>0,$$ 
hence, by \equ{goes to infinity}, $\sup_{Y\in B}\|\rho\big(g_{\vt_n} \tau(Y)\big) v_n\| 
\to \infty$, a contradiction. 
\qed

\section{Theorem \ref{thm: friendly general} and quantitative nondivergence}
\name{est}


\subsection{A sufficient condition}\name{cond} 
The second goal of this paper is to 
show that sets $\DI_\vre(\mathcal{T})$ are null with respect to
 certain
 measures $\mu$ on $\mr$ other than Lebesgue. We will use Proposition  
\ref{prop: dynamical interpretation} to 
formulate a condition sufficient for having $\mu\big(\DI_\vre(\mathcal{T})\big) = 0$
for fixed $\vre > 0$ and all unbounded $\TT\subset \fa^+
$.
Similarly to the setup of Theorem \ref{thm: friendly general},
we will consider measures $\mu$ of the form 
 $F_*\nu$, where $\nu$ is a measure on $\R^d$  and
 $F$ a map from $\R^d$ to $\mr$.

\ignore{
Now suppose that a subset 
$\mathcal{T}_0$ of $\fa^+
$, $\vre > 0$  and a 
measure $\mu$ on $\mr$ are fixed. 
We are going to apply 
the above proposition 
to derive a 
condition guaranteeing that  $\mu\big(\DI_\vre(\mathcal{T})\big) = 0$ 
for any  unbounded subset of $\mathcal{T}_0$.
Slighty more generally, we will be considering
pushforwards $F_*\nu$ of 
measures $\nu$ on $\R^d$ by 
maps $F:\R^d\to\mr$.}


\begin{prop}
\name{prop: density} Let   a 
 measure  $\nu$  on
$\R^d$,   an open subset  $U$ of $\R^d$, 
  a map $F: U\to\mr
$,
and 
$\,0 < \vre,c < 1$ be given. 
\commargin{shortened the statement, reducing to just one part}
Suppose that 
 for any ball $B\subset U$ there exists $s
> 0$ such that the inequality 
\eq{measest}{
\nu\left(
\big\{\x\in B : g_{\vt} \bar{\tau}\big(F(\x)\big)
\notin K_\vre\big\} \right) \le c \nu(B)}
holds for any $\,\vt\in\fa^+$ with $\|\vt\|
 \ge s$.
Then $F_*\nu\big(\DI_\vre(\mathcal{T})\big) = 0$ for any 
unbounded $\,\mathcal{T}\subset\fa^+
$.
\ignore{
the following are equivalent:

\begin{itemize}
\item[(i)] $F_*\nu\big(\DI_\vre(\mathcal{T})\big) = 0$ for any  unbounded $\mathcal{T}\subset\mathcal{T}_0$;
\item[(ii)] there exists $\,c < 1$ with the following property: $\nu$-almost every point of $U$ has 
a neighborhood $U'$ such that
 for any ball $B\subset U'$ one can find $t_0 = t_0(B) > 0$ with 
\eq{measest}{
\vt\in\mathcal{T}_0,\ \|\vt\| > t_0\Rightarrow
\nu\left(
\big\{\x\in B : g_{\vt} \bar{\tau}\big(F(\x)\big)
\notin K_\vre\big\} \right) \le c \nu(B)\,;}
\item[(iii)] for any $0 < c < 1$ and any ball  $B\subset U$ there exists $t_0$ 
such that \equ{measest} holds.
\end{itemize}}
\end{prop}

\begin{proof} 
\commargin{shortened the proof as well}
Since $\TT$ is unbounded, it follows from   the assumption of the proposition that
 for any ball $B\subset U$ 
 and 
 any positive $t$ one has 
\commargin{replaced $Y$ by $F(\x)$}
$$\nu\Big(\bigcap_{\vt\in\mathcal{T} ,\, \|\vt\| \ge t}\big\{\x\in B : g_{\vt} \bar{\tau}\big(F(\x)\big)
\notin K_\vre\big\} \Big) \le c \nu(B)\,.$$
Therefore, by  \equ{liminf}, $\nu\big(\{\x\in B : F(\x)\in \DI_\vre(\mathcal{T})\}\big) 
\le c \nu(B)$.
In view of a density theorem for Radon measures on Euclidean spaces
\cite[Corollary 2.14]{Mattila}, this forces $F^{-1}\big(\DI_\vre(\mathcal{T})\big)$
to have $\nu$-measure zero. 
\end{proof}

\subsection{A \qn\ estimate}\name{quant} 
\ignore{\combarak{See if you like this lead-in to quantitative
nondivergence. As you wrote in your e-mail there is some overlap
between \S2 and \S3 when it comes to non-divergence. I have chosen not
to emphasize it. This is not just laziness, the sections are now
completely independent and should be kept that way}.} 
The proof \commargin{removed ``both proofs"}  of Theorem \ref{thm: unif distr} 
given in \S\ref{dynamical} relies on Corollary \ref{cor: recurrence}, which is a
quantitative nondivergence estimate for translates of unipotent
trajectories. Estimates of this kind  have their origin in the proof by
Margulis \cite{Mar:non-div} that orbits of unipotent flows do not diverge, see
\cite{handbook} for a historical account. During the last decade, starting from the paper 
\cite{KM}, these techniques were transformed into a powerful 
method
yielding measure estimates as in \equ{measest} for a certain  
broad class of measures $\nu$ and maps $F$. To introduce these we need
to elaborate on the definitions mentioned in the introduction.

If $B = B(\x,r)$ is a ball in $\R^d$ and $c > 0$, $cB$ will denote the 
ball $B(\x,cr)$. 
A measure $\nu$ on
$\R^d$ is said to be {\sl $D$-Federer\/} on 
an open 
$U\subset \R^d$ if
for all balls $B$ centered at $\supp\,\nu
$ with $3B\subset U$ one has
${\nu(3B)}/{\nu(B)} \le D$.

If $\nu$ is a measure on $\R^d$, $B$ a subset of $\R^d$ with $\nu(B)>0$, and
$f$  a real-valued measurable function on
$B$, we let
$$\Vert f \Vert_{\nu, B} \stackrel{\mathrm{def}}{=}\sup_{\x \in B\, \cap\,
\supp\, \nu} |f(\x)|\,.
$$
Given $C,\alpha > 0$, open $U \subset \R^d$
and a measure $\nu$ on
$\R^d$, say that $f:U\to \R$ is
         {\sl $(C,\alpha)$-good on $U$ with respect to
$\nu$\/}
         if for any ball $B \subset U$ centered in $\supp\,\nu$
and any
$\vre > 0$ one has
\eq{def-good}{\nu\big(\{y\in B : |f(y)| < \vre\}\big) \le C
\left(\frac{\varepsilon}{\Vert f\Vert_{\nu, B}}\right)^\alpha{\nu(B)}\,.}

We need to introduce some more notation in order to state a theorem from \cite{friendly}. 
Let
$$
\mathcal{W}\df\text{ the set of proper nonzero rational
subspaces of } \R^{k}\,.
$$
From here until the end of this section, we let $\Vert \cdot \Vert$ stand
for  the  Euclidean norm on $\R^{k}$, induced by the standard inner product
$\langle \cdot,\cdot \rangle$, which we extend from
$\R^{k}$ to its exterior algebra. 
For $V
\in \mathcal{W}$ and $g \in G$, let
$$
\ell_V(g) \stackrel{\mathrm{def}}{=} \Vert g (\vv_1 \wedge \cdots \wedge
\vv_j ) \Vert\,,
$$
where $\{\vv_1, \ldots, \vv_j\}$ is a generating set for $\Z^{k} \cap
V$;
note
that $
\ell_V(g)$ does not depend on the choice of  $\{\vv_i\}$.

\begin{thm}[\cite{friendly}, Theorem 4.3]
\name{thm: friendly nondivergence}
Given $d,k\in\N$ and 
positive constants $C,D,\alpha$,
         there exists $C_1 = C_1(d,k,C,\alpha,D) > 0$
with the following property.
Suppose  a measure  $\nu$  on $\R^d$ is  $D$-Federer on 
a ball $\til B$ centered at  $\supp\,\nu$, $0 < \rho \le 1$,
and $h$ is a continuous map $\til B
\to G$ 
such that
for each $V \in \mathcal{W}$,
\begin{itemize}
\item[(i)]
the function
$ \ell_V\compose {h}$  is $(C,\alpha)$-good on $\til B
$
with respect to
$\nu$,
           \end{itemize}
and
\begin{itemize} \item[(ii)]
\label{item: attain rho}
$\Vert \ell_V \compose {h} \Vert_{\nu,B} \geq \rho$, where $B = 3^{-(k-1)}\til B$.
          \end{itemize}
Then for any $\,0<
\varepsilon
\leq \rho$,
$$
{\nu\big(\big\{\x \in B: \pi\big({h}(\x)\big) \notin K_{\varepsilon}
\big\}\big)}\le C_1
\left(\frac{\vre}{\rho} \right)^{\alpha}{\nu(B)} \,.
$$
\end{thm}

\subsection{Checking (i) and (ii)}\name{checking} 
At this point we restrict ourselves to 
the setup of Theorem \ref{thm: friendly general}, that is
consider
measures on $\R^n\cong M_{1,n}$ of the form  $\vf_*\nu$, where $\nu$ is
a measure on $\R^d$ and $\vf= (f_1,\dots,f_n)$ is a map from an open 
$U\subset\R^d$ with $\nu(U) > 0$ to $\R^n$. 
In order to combine  Proposition \ref{prop: density} with
Theorem \ref{thm: friendly nondivergence}, one needs to work with functions
$ \ell_V\compose {h_\vt}$ for each $V \in \mathcal{W}$, where 
\eq{def h}{h_\vt \df  g_\vt\circ\tau\circ\vf\,,}
 and  find conditions sufficient for the validity of (i) and (ii) of 
Theorem \ref{thm: friendly nondivergence} for  large enough $\vt\in\TT$.



The explicit computation that is reproduced below first appeared in \cite{KM}.
 Let $\ve_0,\ve_1,\dots,\ve_n$ be the standard 
 basis  of $\R^{n+1}$, and for
\eq{def I}{I =
\{i_1,\dots,i_{j}\}\subset \{0,\dots,n\},\quad i_1 < i_2 < \dots < i_{j}\,,}  let 
$\ve_{I} \df 
\ve_{i_1}\wedge\dots\wedge \ve_{i_{j}}$;
then $\{\ve_{I}\mid \# I = {j}\}$ is an orthonormal basis of
$\bigwedge^{j}(\R^{n+1})$. Similarly, it will be convenient to put  $\vt = 
(t_0,t_1,\dots,t_{n})\in\fa^+$ where   
\eq{sumequal1}{
\fa^+  = \left\{ 
(t_0,t_1,\dots,t_{n})\in \R^{n+1} : t_i > 0,\ 
t_0 =  \sum_{i = 1}^{n} t_i \right\}
\,.
} 
Then one immediately  sees that for any $I$ as in \equ{def I}, 
\eq{eigenbasis}{\ve_{I}\text{ is an
eigenvector for }g_\vt\text{ with eigenvalue }e^{t_I}\,,} where 
\begin{equation*}
t_I\df  \begin{cases} & t_0 - \sum_{i\in I\nz} t_i \hskip .7in \text{ if }0\in I\\
&- \sum_{i\in I} t_i\ \hskip 1in \text{ otherwise.}\end{cases}
\end{equation*}
We remark that in view of \equ{sumequal1}, $e^{t_I}$ is not less than $1$ 
 for any $\vt \in\fa^+$ and $0\in I$.
Moreover, let $\ell$ be such that $t_\ell = \max_{i = 1,\dots,n}t_i$. 
Then
\eq{big eigenvalue}{{t_I}\ge t_\ell \ge \|\vt\|/n \quad
 \forall\,I
 \text{ containing }0\text{ and not containing }\ell\,.}

Since the action of
$\tau(\y)$, 
where $\y\in\R^n$,
leaves $\ve_0$ invariant and sends $\ve_{i}$, $i > 0$,  to $\ve_{i} + y_i
\ve_0$, one can write\footnote{The choice of 
  $+$ or $-$ in \equ{computation ei} depends on the parity of the number of elements of $I$ 
less than $i$ and is not important for our purposes. See however \cite{dima exponents} for a
more precise  computation.}
\eq{computation ei}{
\tau(\y)\ve_{I} =  \begin{cases} &\ve_{I} \hskip 1.58in \text{ if }0\in I\\
&\ve_{I} + \sum_{i\in I} \pm y_i\,
\ve_{I \cup \{0\}\ssm\{i\}}\ \text{ otherwise.}\end{cases}
}
Now  take  $V
\in \mathcal{W}$, choose  a generating set $\{\vv_1, \ldots, \vv_j\}$ for $\Z^{k} \cap
V$, and expand $\vw \df \vv_1 \wedge \cdots \wedge
\vv_j$ with respect to the above basis by writing  $\vw  = 
 \sum_{I\subset \{0,\dots,n\},\,\# I = j}w_{I}\ve_{I}\in \bigwedge^j(\Z^{k})\nz$. Then one has
\eq{computation w}{
\tau(\y)\vw  = \sum_{0\notin I}
w_{I}\ve_{I} + \sum_{0\in I}\left( w_{I} + \sum_{i\notin I} \pm w_{I\cup\{i\}\ssm\{0\}}y_i \right) 
\ve_{I}\,.
}
\ignore
{hence, for any $\vt$
and any $I$, the projection 
of 
$g_\vt\tau
(\y)\vw$ onto $\ve_I$ is equal to
\eq{computation h}{
\langle g_\vt\tau
(\y)\vw,\ve_I\rangle = 
\begin{cases} & e^{t_I}w_{I} \hskip 1.8in \text{ if }0\notin I\\
&e^{t_I}\left( w_{I} + \sum_{i\notin I} 
\pm w_{I\cup\{i\}\ssm\{0\}}y_i \right) \ \text{ otherwise.}\end{cases}
}}
Here is an immediate implication of the above formula:

\begin{lem}\name{lem: ext power affine} For any $\vw$ 
and $\vt$, 
the map $\y\mapsto g_\vt\tau
(\y)\vw$ is affine; in other words,  for any $I\subset \{0,\dots,n\}$  the projection 
of 
$g_\vt\tau
(\y)\vw$ onto $\ve_I$ has the form
\eq{lincomb}{\langle g_\vt\tau
(\y)\vw,\ve_I\rangle =  c_0 + \sum_{i = 1}^nc_iy_i}
 for some $c_0,c_1,\dots,c_n \in \R$.
\end{lem}


It is also clear from \equ{computation w} that $\tau(\cdot)\vw  \equiv \vw$ if the subspace $V$ represented by $\vw$
contains $\ve_0$ (in other words, if $w_I = 0$ whenever $0\notin I$). In  this case
for any $I\ni 0$  all the coefficients $c_i$ in \equ{lincomb} with $i \ge 1$ are equal to $0$, and $|c_0| = e^{t_I}|w_I| \ge 1$ as long as $w_I \ne 0$.

Let us now consider the complementary case. 

\begin{lem}\name{lem: ext power big} Suppose that $\vw \in \bigwedge^j(\Z^{k})$,  $1\le j \le n$,
is such that  
 $w_J \ne 0$ for some $J\subset \{1,\dots,n\}$.
Then for any  $\vt\in\fa^+$
 there  exists  $I$ 
such that the absolute value of one of
the coefficients $c_i$ in \equ{lincomb} is at least $e^{\| \vt\|/n}$.
\end{lem}

\begin{proof} Let $\ell$ be such that $t_\ell = \max_{i = 1,\dots,n}t_i$. If $J$ as above contains
$\ell$, take  $I = J\cup\{0\}\ssm\{\ell\}$. Then $J = I\cup\{\ell\}\ssm\{0\}$, hence, by \equ{computation w}
and \equ{eigenbasis}, the coefficient $c_\ell$ in the expression  \equ{lincomb} 
for this $I$ is equal to $\pm  e^{t_I} w_J$. The claim then follows from \equ{big eigenvalue},
since $|w_J| \ge 1$ and $I$ contains 0 and does not contain $\ell$. 

If $J$ does not contain $\ell$, the argument is similar: choose any $i\in J$ and take $I = J\cup\{0\}\ssm\{i\}$. Then $J = I\cup\{\i\}\ssm\{0\}$, hence the coefficient $c_i$ in   \equ{lincomb} 
for this $I$ is equal to $\pm  e^{t_I} w_J$. As before, $I$ contains 0 and does not contain $\ell$,
so  \equ{big eigenvalue} applies again.
\end{proof}

\subsection{Proof of Theorem \ref{thm: friendly general}}\name{applying} 
The next theorem generalizes \cite[Theorem 5.4]{KM}:

\begin{thm}\name{thm: nondiv} For any\,
$d,n\in\N$ and  any 
$\,C,\alpha,D > 0$ there exists
$C_2 =  C_2(d,n,C,\alpha,D)$ with the following property.
Suppose  a measure  $\nu$  on $\R^d$ is  $D$-Federer on 
a ball $\til B$ centered at  $\supp\,\nu$ and 
$\vf:\tilde B\to\R^n$  is continuous. Assume that:
\begin{itemize}
\item[(1)] 
any linear 
combination of
$1,f_1,\dots,f_n$  is
$(C,\alpha)$-good on
$\tilde B$ with respect to $ \nu$;
\item[(2)] 
the
restrictions of $ 1,f_1,\dots,f_{n}$ to $B\,\cap \,\supp\,\nu$,  where $B = 3^{-n}\til B$,
are linearly independent over $\R$.
           \end{itemize}
Then there exists $s  > 0$ such that 
 for any $\,\vt\in\fa^+$ with $
\|vt\| \ge s$ and any $\,\vre < 1$, one has 
$$
{\nu\big(\big\{\x \in B: g_\vt\bar\tau\big(\vf(\x)\big) \notin K_{\varepsilon}
\big\}\big)}\le C_2
\vre^{\alpha}{\nu(B)} \,.
$$
\end{thm}

\ignore{\commargin{I think in some form we should relate Lemma \ref{lem: ext power}(b) with 
what we did in \S\ref{drift}, and Theorem \ref{thm: nondiv} with Corollary
\ref{cor: recurrence}. \combarak{I disagree, see comment at the beginning of this
section}}}
\begin{proof}
We will apply Theorem \ref{thm: friendly nondivergence} with $h = {h_\vt}$
as in \equ{def h}. Take  $V
\in \mathcal{W}$ and, as before,  represent it by $\vw = \vv_1 \wedge \cdots \wedge
\vv_j$, where  $\{\vv_1, \ldots, \vv_j\}$ is a generating set for $\Z^{k} \cap
V$. From Lemma \ref{lem: ext power affine} and assumption (1) above it follows that
for any $\vt$, each coordinate 
of $h_\vt(\cdot)
\vw$ 
is $(C,\alpha)$-good on
$\tilde B$ with respect to $ \nu$. Hence the same, with $C$ replaced by 
${n+1\choose j}^{\alpha/2}  C$, can be said about
 $\Vert h_\vt(\cdot)
\vw\Vert = \ell_V\circ h_\vt$, see \cite[Lemma 4.1]{friendly}.
This verifies condition (i) of Theorem \ref{thm: friendly nondivergence}.

Now observe that  assumption (2) implies the existence of $\delta > 0$ (depending on $B$)
such that
$\|c_0 + \sum_{i = 1}^nc_if_i \|_{\nu,B} \ge \delta$ for any $c_0,c_1,\dots,c_n 
$
with $\max |c_i| \ge 1$.  Using Lemma \ref{lem: ext power big} and the remark preceding it,
 we conclude that either
$\|\ell_V\circ h_\vt \|_{\nu,B} \ge 1$ (in the case $\ve_0\in V$) or $\|\ell_V\circ h_\vt \|_{\nu,B} \ge\delta e^{\| \vt\|/n}$ (in the complementary case).
So  condition (ii) of Theorem \ref{thm: friendly nondivergence} holds with $\rho = 1$
whenever $\| \vt\|$ is at least $s \df  - n\log \delta$.
\end{proof}

We can now proceed with the 

\begin{proof}[Proof of Theorem \ref{thm: friendly general}]
It suffices to show that for $\nu$-a.e.\ 
$\x$  there exists  a ball $B$ centered at $\x$ such that 
\eq{measure 0}{\nu\big(\{\x\in B : \vf(\x) \in \DI_\vre(\mathcal{T})\}\big) = 0\,.} 
Since $(\vf,\nu)$ is 
$(C,\alpha)$-good and
nonplanar, for $\nu$-a.e.\ 
$\x$ one can choose $B$ centered at $\x$ such that  $(\vf,\nu)$ is nonplanar on $B$ 
and \cag\ on $\tilde B = 3^nB\subset U$, which implies that conditions (1) and 
(2) of Theorem \ref{thm: nondiv}
are satisfied. 
Then, combining  Theorem \ref{thm: nondiv} with 
Proposition \ref{prop: density}, one concludes that \equ{measure 0} holds whenever
$C_2\vre^\alpha < 1$. \end{proof}

We remark that the constant $C_1$ from 
Theorem \ref{thm: friendly nondivergence}, and hence $C_2$
from  Theorem \ref{thm: nondiv} and $\vre_0$
from  Theorem \ref{thm: friendly general}, can be explicitly
extimated in terms of the input data of those theorems, 
see \cite{KM, BKM, KT, dima  exponents}.  However we chose not to bother
the reader with
explicit computations, the reason being that in the special cases previously considered in the
 literature
our method produces much weaker estimates. More on that in the next section.

\pagebreak

\ignore{The goal of this section is to show that the set of matrices for which $(\vr,\vs)$-DT 
{\sl can be $\delta$-improved  
 along\/} a subsequence $\mathcal{T}$ is null with respect to some measures 
other than Lebesgue. 

 sets $\DI_\delta(\vr,\vs;\mathcal{T})$ 
has Lebesgue 
measure zero
Using an elementary argument which can be found e.g.\ in  \cite[Chapter V, \S 7]{Cassels}
and dates back to  Khintchine, \commargin{(I think)}
 one can show that 
for any $\vr,\vs$ and $\mathcal{T}$
as above,  
$\DI_\delta(\vr,\vs;\mathcal{T})$ has Lebesgue 
measure zero as long as $\delta < 2^{-(m+n/2)}
$. 

At this point it is worthwhile to point out 
that the definition \equ{keps} of $K_\vre$ depends on the choice of the norm on 
$\R^{n+1}$. In this section it will be convenient to use Euclidean norm instead
of the supremum norm. The only way it will affect the proof of Theorem 
\ref{thm: friendly general} is through...

}

\section{Examples and applications}
\name{const}

\subsection{Polynomial maps}\name{polynom} A model example 
of functions which are \cag\ with respect to Lebesgue measure is given 
by polynomials: it is shown in \cite[Lemma 3.2]{KM} that any polynomial of
degree $\ell$ is \cag\ on $\R$ with respect to $\lambda$, where $\alpha = \frac1\ell$ and 
$C$ depends only 
on $\ell$. The same can be said about polynomials in $d\ge 1$ variables, with 
$\alpha = \frac1{d\ell}$  and $C$ depending on
 $\ell$ and $d$. Obviously $\lambda$ is $3^d$-Federer
on $\R^d$. Thus, as a corollary of Theorem \ref{thm: friendly general}, we obtain
the existence of $\vre_1 = \vre_1(n,d,\ell)$ such that 
whenever $\vre < \vre_1$, \equ{result} holds for
$\nu = \lambda$ and any polynomial map 
$\vf = (f_1,\dots,f_n)$ of degree $\ell$  in $d$ variables such that
$ 1,f_1,\dots,f_{n}$
are linearly independent over $\R$. This in particular applies to 
$\vf(x) = (x,\dots,x^n)$, a generalization of the setup of Theorem \ref{thm: ds2}
considered by Baker in \cite{Baker JLMS} for $n = 3$ and then by 
Bugeaud for an arbitrary $n$. Note that it is proved 
in \cite{Bugeaud} that $(x,\dots,x^n)$ is almost surely 
not in $\DI_\vre$ for $\vre < 1/8$. Our method, in comparison, shows that
$(x,\dots,x^n)$ 
is almost surely
not in $\DI_\vre(\TT)$ for any $\TT$ drifting away from walls
and  $\vre < 1/n^n (n+1)^{2} 2^{n^2 + n}$.

Improving these results to any $\vre<1$ is a natural and challenging
problem. Recently \cite{Nimish new} the following was obtained: 
\begin{thm}[N.A.\ Shah]
\name{thm: Nimish}
Let $\vf: \R \to M_{1,2} \cong \R^2$ be a nonplanar 
 polynomial
map. For an interval $B 
\subset \R$ and for ${z} \in \SL_3(\R)/\SL_3(\Z)$, let $\nu_{{z},B}$ be
the natural probability measure on $\tau\big(\vf(B)\big)z$, that is the pushforward of
$\frac1{\lambda(B)} \lambda|_B$ via the map $x \mapsto \tau \circ \vf(x){z}.$ Then for any
unbounded sequence $\vt_n \in \fr$, the sequence of translated measures
$(g_{\vt_n})_* \nu_{z,B}$ weak-$*$ converges to \rm{vol}.

\end{thm}

The proof follows a similar strategy as our proof of Theorem \ref{thm: unif
distr}, but is
considerably more difficult. 
Repeating the argument of \S 2.2, one obtains:
\begin{cor}
For any $\vre <1$ and any nonplanar polynomial curve $\vf: \R \to
M_{1,2}$,  $\vf(x)$ does not belong to
$\DI_\vre$ for $\lambda$-a.e.\ $x \in \R$. 
\end{cor}

Note also that Bugeaud's result mentioned above can be rephrased in terms of 
small values of integer polynomials at almost all real $x$. 
Theorem 1.5 produces a similar result for polynomials of $d$ variables. 
Namely, for some $c = c(d,n)$ and for $\lambda$-almost every
$\x\in\R^d$ there are infinitely many integers $N$ for which there are no polynomials 
$P\in\Z[X_1,\dots,X_d]$  of degree at most $n$ and height  less than $ N$  with $
|P(\x)| < c N^{- m}
$,
where $m$ is the dimension of the space of non-constant
polynomials in $d$ variables  of degree at most $n$. 

\subsection{Nondegenerate maps}\name{nondeg}
Here is another situation in which \cag\ functions arise. The following lemma
is a strengthening of  \cite[Lemma 3.3]{KM}: 

\begin{lem}\name{lem: nondeg good} For any $\,d\in\N$ there exists
$C_{d} > 0$ with the following property.
Let $B$ be a 
cube in $\R^d$ (product of intervals of the same length), 
and let 
$f\in C^\ell(B)$, $\ell\in\N$, be such 
that for some positive constants $a_1,\dots, a_d$ and $A_1,\dots, A_d$ one has 
\eq{bounds}{
a_i\le |\partial_i^{\ell}f(\x)|\le A_i\quad \forall\,\x\in B,\ i = 1,\dots,d\,.
}
Then 
$$
\lambda\big(\{\x\in B: |f(\x)| < \vre\}\big) \le \ell C_{d}
\max_i \left(\frac{A_i}{a_i}\right)^{1/\ell} \left(\frac\vre{\|f\|_{\lambda, B}}\right)^{1/d\ell}
\lambda(B)\,.
$$
\end{lem}

We remark that \cite[Lemma 3.3]{KM} instead of \equ{bounds} assumed
$$
a\le |\partial_i^{\ell}f(\x)|\quad \forall\,\x\in B,\ i = 1,\dots,d
$$
and
$$
|\partial^{\beta}f(\x)|\le A\quad \forall\,\x\in B \ \forall \text{ multiindex }\beta 
\text{ with }|\beta|\le \ell\,,
$$
and produced the same conclusion as the above lemma, 
with $A_i = A$
and $a_i = a$ for all $i$.

\begin{proof}[Sketch of Proof] The case $d = 1$ can be proved by a verbatim repetition
of the  argument from  \cite{KM} -- it is easy to verify that
a bound on just the top derivative  is enough for the proof.
The general case then follows using \cite[Corollary  2.3]{KT}.
\ignore{, which 
in particular states that whenever $f$ is a function
on $U_1\times\dots\times U_d$ ($U_i\subset \R$ open) 
such that for any $j = 1,\dots,d$ and
any $x_i \in U_i$ with $i \ne j$, the function
$$
y\mapsto f(x_1,\dots,x_{j-1},y,x_{j+1},\dots,x_d)
$$
is $(C,\alpha)$-good on $U_j$, one knows that 
$$f\text{ is $(d C, \alpha/d)$-good  on }
U_1\times\dots\times U_d
$$
 with respect to the product metric. }
\end{proof}

Recall that a map $\vf$ from  $U\subset\R^d$ to $\R^n$  is called 
{\em $\ell$-nondegenerate at\/}  $\x\in U$ if   partial derivatives of
$\vf$ at $\x$
up to order $\ell$ span $\R^n$, and {\em $\ell$-nondegenerate\/}
if it is $\ell$-nondegenerate at $\lambda$-a.e.\ $\x\in U$. Arguing as in
the proof of \cite[Proposition 3.4]{KM}, from the above lemma one deduces

\begin{prop}\name{prop: nondeg good} For any $\,d,\ell\in\N$ there exists
$C_{d,\ell} > 0$ with the following property.
Let $n\in\N$ and let $\vf = (f_1,\dots,f_n):U\to \R^n$ be $\ell$-nondegenerate at 
$\x\in U\subset \R^d$. 
Then for any $C > C_{d,\ell}$ there exists a neighborhood 
$V\subset U$ of $\x$ 
such that any linear combination of $ 1,f_1,\dots,f_n$ is $(C,1/d\ell)$-good on 
$V$. 
\end{prop}

This implies that if $\vf:U\to \R^n$ is $\ell$-nondegenerate, then 
the pair $(\vf,\lambda)$ is 
$(C,1/d\ell)$-good for any $C > C_{d,\ell}$. Also,  nonplanarity is
clearly  an 
immediate consequence of nondegeneracy. Thus, by  Theorem \ref{thm: friendly general},
there exists 
$\vre_2 = \vre_2(n,d,\ell)$ such that 
whenever $\vre < \vre_2$,  \equ{result} holds for
$\nu = \lambda$ and any  $\ell$-nondegenerate
$\vf:U\to \R^n$, $U\subset \R^d$.  This was previously established in the case
$d = 1$, $l = n = 2$ in \cite{Baker MPCPS}, with an additional 
assumption that $\vf$ be $C^3$ rather than $C^2$.  Also, 
M.\ Dodson, B.\ Rynne, and J.\ Vickers
considered $C^3$ 
submanifolds of $\R^n$ 
with `two-dimensional definite curvature almost
everywhere', a condition which implies $2$-nondegeneracy 
(and requires the dimension of  the manifold   to be at least $2$).
It is proved 
in \cite{DRV1} that 
almost every point on such a manifold 
is 
not in $\DI_\vre$ for $\vre < 2^{-\frac n {n+1}}$.  We also remark that 
the result of  \cite{Nimish new} extends to nondegenerate analytic curves in $\R^2$.

\subsection{Friendly measures}\name{fr} The class of friendly measures
was introduced in \cite{friendly}, the word `friendly' being
an approximate abbreviation of  `Federer, nonplanar and decaying'. Using
the terminology of the present paper, we can define this class as follows:
a measure $\mu$ on $\R^n$ is  {\sl friendly\/} if for $\mu$-a.e.\ 
$\x\in\R^n$  there exist a neighborhood $U$ of $\x$ and  $D, C, \alpha > 0$ 
such that $\mu$
is
 $D$-Federer  on $U$, and $(\Id,\mu)$ is  both  \cag\ and 
nonplanar\footnote{We remark that the nonplanarity of $(\Id,\mu)$ as defined 
in this paper
 is a condition weaker than the nonplanarity of $\mu$ as defined 
in \cite{friendly}.} on $U$.
%
In order to apply  Theorem \ref{thm: friendly general} we would 
like to use somewhat more uniform version:
given $C, \alpha, D > 0$, define
 $\mu$   to be   {\sl 
$(D,C,\alpha)$-friendly\/} if  for $\mu$-a.e.\ 
$\x\in\R^n$  there exists a neighborhood $U$ of $\x$  
such that $\mu$
is
 $D$-Federer  on $U$  and $(\Id,\mu)$ is both  \cag\ and nonplanar on $U$.
In view of Theorem \ref{thm: friendly general}, 
almost all points with respect to those measures are  not in $\DI_\vre(\TT)$ 
for any $\TT$ drifting away from walls
and small enough $\vre$, where $\vre$ depends only on $C,
D, \alpha$.

As discussed in the previous subsections, smooth measures
on nondegenerate submanifolds of $\R^n$ satisfy the above properties. 
Furthermore, the class of friendly measures is rather large; 
many examples are described in \cite{friendly, bad,
Urbanski, Urbanski - d=1, Urbanski-Str}. 
A notable class of examples is given by limit measures of finite 
irreducible systems of contracting similarities \cite[\S8]{friendly}
(or, more generally, self-conformal contractions, \cite{Urbanski}) 
of $\R^n$ with the
open set condition.
These measures were shown to be $(D,C,\alpha)$-friendly for some 
$D,C,\alpha$, thus satisfy the conclusions of Theorem \ref{thm: friendly general}.

\subsection{Improving DT along non-drifting $\TT$}\name{nodrift} 

Comparing Theorem \ref{thm: dsw} with Theorem \ref{thm: friendly general}, one 
sees that the former has more restrictive assumptions, namely  $\TT$ has
to drift away from walls as opposed to just be unbounded. This is not an accident: 
the drift condition is in fact necessary for the main technical tools of the proof,
that is,  equidistribution results of \S\ref{dynamical}. 

To see this, for simplicity 
let us restrict ourselves to the case $m = 2$, $n = 1$; 
the argument for the general case is similar.
Suppose that $\TT = \{ \T^{(\ell)}  = (t_1^{(\ell)},t_2^{(\ell)},t_3^{(\ell)} ): \ell\in\N\}\subset \fa^+$ is unbounded but
does not satisfy \equ{drift}; that is, either  
$\{t_1^{(\ell)}\}$ or  $\{t_2^{(\ell)}\}$ is bounded.
Without loss of generality, and passing to a subsequence, we can assume that  $t_1^{(\ell)}$ 
is convergent as $\ell\to\infty$; that is, for any $\ell$
we can write  $\T^{(\ell)} = \vs^{(\ell)}  + \vu^{(\ell)} $ where   
\eq{conditionons}{\vs^{(\ell)} = (0,s^{(\ell)},s^{(\ell)}),\quad 
s^{(\ell)} \to\infty\text{ as }\ell\to\infty\,,}
 and $g_{\vu^{(\ell)}}\to g_0\in \SL_3(\R)$ as $\ell\to\infty$.
Note that for any $\ell$ and any $Y\in M_{2,1}$, $g_{s^{(\ell)}}\tau(Y)$ belongs to the
group 
\eq{def H}{H \df \left\{\begin{pmatrix} 1 & \x\\ 0 & g\end{pmatrix} : g\in\SL_2(\R),\ \x\in\R^2\right\}\,,}
a semi-direct product of $\SL_2(\R)$ and $\R^2$. 
Thus for any $Y$, the trajectory 
$\{ g_{\T}\tau(Y)\Z^3: \T\in\TT\}$ must approach $g_0H\Z^3$, 
which is a proper submanifold  inside the space of lattices in $\R^3$
(in fact, $H\Z^3$ is the set of lattices in $\R^3$ containing $\ve_1 = (1,0,0)$
as a primitive vector); therefore the translates $g_\T\tau(B) \Z^3$, where $B$ is 
any subset of $M_{2,1}$,  do not become equidistributed.
We conclude that it is not possible to prove the analogue of Theorem \ref{thm: unif distr}
with $\lfloor \T\rfloor$ replaced by $\|\T\|$.

 \medskip
Note that this apriori does not rule out proving Theorem  \ref{thm: dsw} with a relaxed assumption on $\TT$:
recall that our goal was to make almost every orbit return to a specific set $K_\vre$, not just 
any nonempty open set.
And indeed, it is easy to show, using Theorem  \ref{thm: dsw}, induction and Fubini's Theorem, that 
the set $\DI_\vre(\TT)$ is $\lambda$-null whenever  $\vre < 1$ and $\TT$ is such that 
for every $i$, $\{t_i : \T \in \TT\}$
 is either unbounded or converges to $0$.

 However, in general Theorem \ref{thm: dsw} is false if
 one just assumes that $\TT$ is unbounded. 
 Here is a simple counterexample, also in 
 the case $m = 2$, $n = 1$. 
Fix  $0 < \vre < 1$ and write  $\T^{(\ell)} = \vs^{(\ell)}  + \vu $, where   $ \vs^{(\ell)} $ is as in \equ{conditionons}
and $\vu = (u,u,0)$. Then for any $\ell$ and any $Y$,  $g_{ \vs^{(\ell)}}\tau(Y)\Z^3$ belongs to 
$H\Z^3$ where $H$ is as in \equ{def H}, that is,  contains $\ve_1$
as a primitive vector. Therefore for any $\T\in\TT$ and any $\s$,   $g_{ \T}\tau(Y)\Z^3$
contains $e^u\ve_1$
as a primitive vector. Now suppose that $1/\vre^2 < e^u < 2\vre$.
Let $B_\vre\subset \R^3$ be given by $\{|x_1| < e^u, |x_2| < \vre, |x_3| < \vre\}$.
It is a convex centrally symmetric domain of volume greater than $ 8$, 
hence, by Minkowski's Lemma,
 it
must contain a nonzero  vector $\vv\in g_{ \T}\tau(Y)\Z^3$. However, since $e^u < 2\vre$,
the sup-norm distance of $\vv$ to either  $e^u\ve_1$ or $-e^u\ve_1$ is less than $ \vre$,
and it is positive since $\pm e^u\ve_1\notin B_\vre$.
This proves that $g_{ \T}\tau(Y)\Z^3$ is always disjoint from $K_\vre$; thus, under those
assumptions on $\TT$ and $\vre$, the set $\DI_\vre(\TT)$ is equal to $M_{2,1}$.

Similar counterexamples exist in any dimension. Still, it seems plausible that Theorem  \ref{thm: dsw} 
will remain valid if  \equ{drift} is replaced by an assumption that $\liminf_{\ell\to\infty}\lfloor \T^{(\ell)} \rfloor$
is large enough. The proof of this requires equidisiribution results for more general homogeneous spaces (specifically, spaces similar to $H\Z^3$ in the above example).

\subsection{Weighted badly approximable systems}\name{ba} We conclude the paper
with another application of Theorem \ref{thm: unif distr}. 
Let $g_t$ be a one-parameter
subgroup of $G$. Suppose a subgroup $H$ 
of $G$ normalized by $g_t$ is such that (a) the conjugation by $g_t$, $t > 0$,
 restricted to $H$ is an expanding
automorphism of $H$, and (b) $g_t$-translates of the leaves $Hx$, $x\in\ggm$,
become equidistributed as $t\to\infty$, with the convergence uniform
as $x$ ranges over compact
subsets of $\ggm$. These conditions were shown in \cite{KlMa} to imply
that for any $x\in\ggm$, the set 
$$\big\{h\in H : \text{
 the trajectory }\{g_thx : t > 0\} \text{ is bounded}\big\}
$$ is thick (that is, has full \hd\ at every point). When $H$ is as in \equ{defn u} 
and $\{g_t: t > 0\}$ is any one-parameter
subsemigroup of $G$ contained in $\exp(\fa^+)$, both (a) and (b) are satisfied, 
the latter being a consequence  of Theorem \ref{thm: unif distr}.

Let us now take $g_t$ of  the form
\begin{equation*}
\label{eq: new defn g_t}
g_t =  \diag(e^{r_1t}, \ldots,
e^{r_{m} t}, e^{-s_1t}, \ldots,
e^{-s_{n}t} )\,,
\end{equation*}
where 
\eq{def rs}{r_i,s_j > 0\quad\text{and}\quad\sum_{i=1}^m r_i = 1 =
 \sum_{j=1}^n s_j\,.}
Then it is known \cite{dima duke} that the trajectory 
 $\{g_t\bar\tau(Y): t > 0\}$ is bounded in $\ggm$
if and only if $Y$ is {\sl $(\vr,\vs)$-\ba\/}, which by definition means 
$$ \inf_{\vp\in\Z^m,\,\vq\in\Z^n\nz} \max_i |Y_i\vq - p_i|^{1/r_i} \cdot
\max_j|q_j|^{1/s_j}    > 0\,.
$$
The components of vectors $\vr,\vs$ should be thought of as weights assigned to
linear forms $Y_i$ and integers $q_j$.
Thus one can obtain a weighted generalization
of W.M.\ Schmidt's theorem \cite{Schmidt} on the thickness of the set
of badly approximable systems of linear forms:

\begin{cor}\name{cor: ba} For any choice of $\vr,\vs$ as in \equ{def rs},
the set of $(\vr,\vs)$-\ba\/ \amr\ is thick. 
\end{cor}

This was previously established by A.\ Pollington and S.\ Velani \cite{PV-bad}
in the case $n = 1$.

\end{document}